# Variational Iteration Method for Image Restoration

Keyvan Yahya, Jafar Biazar, Hossein Azari, Pouyan Rafiei Fard

*Abstract*— The famous Perona-Malik (P-M) equation which was at first introduced for image restoration has been solved via various numerical methods. In this paper we will solve it for the first time via applying a new numerical method called the Variational Iteration Method (VIM) and the correspondent approximated solutions will be obtained for the P-M equation with regards to relevant error analysis. Through implementation of our algorithm we will access some effective results which are deserved to be considered as worthy as the other solutions issued by the other methods.

## I. Introduction

SINCE the 90's, PDE based methods have been welcomed widely by vision researches in image restoration and the other applications like Image Segmentation, Motion Analysis and Image Classification and etc. These methods are divided into three categories namely the second order PDEs, the fourth order PDEs and the complex diffusion [1]-[4]. In 1990, Pietro Perona and Jitendra Malik in their common paper [5] for the first time introduced a parabolic partial differential equation which was called thereafter "Perona-Malik" and we abbreviated it by "PM". This equation could be posed rapidly as a new mathematical tool with the great capabilities to image restoring by convolving the original image with a Gaussian kernel. The original form of the equation in [5] have been shown as:

$$\frac{\partial u(x,y,t)}{\partial t} = \nabla \cdot \big(c(x,y,t)\nabla u(x,y,t)\big)$$
$$u(x,y,0) = u_0(x,y)$$
(1)

Where $c(x,y,t)$ is the diffusion factor which were suggested in the two forms of $c(x,y,t) = \frac{u}{(1+\frac{|\nabla u|^2}{k^2})}$

Keyvan Yahya is now an independent researcher in the field of computational neuroscience and applied mathematics (corresponding author to provide e-mail: keyvan.yahya@aol.com).
Jafar Biazar, is now with the Mathematics Department, University of Guilan, P.O.Box 1914, P.C. 4193833697, Rasht, Iran (e-mail: biazar@guilan.ac.ir.).
Hossein Azari is with Faculty of Mathematical Science, University of Shahid Beheshti P.O.Box 19835-389, Tehran, Iran (e-mail: h_azari@sbu.ac.ir).
Pouyan Rafiei Fard is now an independent researcher in the field of Neuromusicology and computational neuroscience. (email: rafieifard@ce.sharif.edu)

and $c(x,y,t) = \exp{(-\frac{|\nabla u|^2}{(k^2)})}$ and also $u_0$ is the original image should be restored.

Some people tried to solve the P-M by various numerical methods such as finite difference method [6], explicit finite volume method [7], total variation method [8] and etc. Each of these methods have been used with regards to some different criterions in empirical implementations for example in edge conservation, speed of algorithm and how optimized the algorithms is. They also showed its considerable efficiency via comparison between former filter methods and these PDEs. Recently, another solution was given by the authors based on homotopy perturbation method [9]. of course, homotopy perturbation method differs to some extent from the other methods we mentioned above, because it is indeed a mesh-free method, like the VIM, while the others are not [10].

Through the studying of a non-linear parabolic equation of the following general form, we search for a new approximated solution:

$$\begin{cases} \frac{\partial u(t,x)}{\partial t} = g(|G_\sigma * \nabla u|)|\nabla u| div \frac{\nabla u}{|\nabla u|}, \\ \frac{\partial u(t,x)}{\partial N} = 0 \\ \frac{\partial u(0,x)}{\partial t} = u_0(x) \end{cases}$$ (2)

where $u(t,x)$ is the solution of this PDE (restored image) we are searching for, $N$ is an outward Normal to domain $\Omega$ and $u(0,x)$ is the original image. The gray-level function $u(t,x)$ is depended on two parameters; the scale parameter denoted by $t$ and the spatial coordinate $x$. $G$ is smoothing Kernel (for instance, an additive Gaussian filter) and $g$ is a non-increasing Lipchitz function.

But to get more conceivable results, many authors manipulated this general form of restoration equation and expanded the P-M to propose different forms of this equation. For example, Lions and Alvarez offered an interesting non-linear form of restoration equation [6] or in [11], the Authors also gave a different parabolic equation which is considerable too or in [12]-[13], the work was focused on a complex form of diffusion equation and so on.

To reach our aim, we appoint to use the non-linear and modified form of P-M equation proposed by [14] because of

its facility and vigor. The equation they proposed can be expressed as the following form:

$$\begin{cases} \frac{\partial u(t,x)}{\partial t} = -F(x, u, \nabla u, \nabla^2 u), \\ \frac{\partial u(t,x)}{\partial N} = 0, \quad in\ [0,t] \times \partial\Omega \\ \frac{\partial u(0,x)}{\partial t} = u_0(x), \quad in\ \Omega \end{cases} \quad (3)$$

where $F(x, u, \nabla u, \nabla^2 u) = -div\ (g(|G_\sigma * \nabla u|)|\nabla u|)$, $\Omega \subset \mathbb{R}^2$ is a rectangular domain, g is a non-increasing Lipchitz function satisfies this condition: $g(0) = 1$, $0 < g(s) \rightarrow 0$ for $s \rightarrow \infty$ and $G_\sigma \in C^\infty(\mathbb{R}^d)$ is a smoothing Kernel with compact support with:

$$\int_{\mathbb{R}^d} G_\sigma(x) = 1\ \text{and}\ G_\sigma(x) \rightarrow \delta_z\ \text{as}\ \sigma \rightarrow 0. \quad (4)$$

$\delta_z$ is the Dirac delta function and because of boundary condition, we can assume that $u_0(x) \in W^{1,2}(\Omega)$ [15].

We are seeking for a solution belongs to a functional space which is the intersection of the Sobolov and the Lebesgue spaces ($L^1(\Omega)$). thus, we can consider $u$ as:

$$u \in W^{1,2}(\Omega) \cap L^1(\Omega), \quad (5)$$

As we mentioned before, most of the methods have been applied in order to find a solution ensue a point by point process of approximation. In fact, the VIM has a more practical essence. To work with it, first of all we must construct *iterative series* in which we impose our PDE.

To start our study about the new method, let us express the explicit form of Perona-Malik equation used to restore the original image. The non-linear form of P-M equation can be expressed as the following:

$$\begin{cases} \frac{\partial u}{\partial t} = g(\|G_\sigma * \nabla u\|)|\nabla u|div\frac{\nabla u}{|\nabla u|}, \\ \frac{\partial u(t,x)}{\partial \eta} = o, \\ \frac{\partial u(0,x)}{\partial t} = u_0(x), \end{cases} \quad (6)$$

where $\eta$ is a normal outward to the domain $\Omega$, and as we know the function $u(t, x)$ is the solution of the nonlinear parabolic equation (restored image) we are seeking for it based on $u(0, x)$(original image). to work with this PDE, we are involved with three specific parameters, the first is the scale parameter denoted by $t$ (along with $x$ in spatial domain) and the second is called "Gaussian Kernel" sometimes we take it as the *additive Gaussian filter* and ultimately $g(.)$ is an non-increasing Lipchitz function. It is the non linear one which enhances edges temporarily, before slowly blurring them out. It may be ill-posed, but can be made well-posed by replacing |r u| by |r uσ| where uσ is the convolution of $u$ by a Gaussian of standard deviation σ. However it is numerically well-posed, because of the implicit regularizing effect of numerical derivatives.

For numerical implementation we must adapt trivially the scheme proposed for the edge preserving linear diffusion but re-computing the diffusivity at each step. Homotopy Method is a Free Mesh numerical instrument which can give the solution. We put our PDE in the algorithm that we will offer and gain a solution which is comparable with the others obtained via other mesh free methods. Also, the existence and uniqueness of the solution gained by the VIM was proved earlier [16]. Of course, many authors see the VIM as a tool for finding an appropriate starting point $x^0$ in Newtonian method which is still a crucial problem.

## II. BASIC IDEA OF VARIATIONAL ITERATION METHOD

Variational iteration method is a powerful device to solve the various kinds of linear and non-linear functional equations. The VIM was developed by He [17]-[19] for solving linear, nonlinear, initial and boundary value problems. But in fact, the first ideas of the method was issued by [20]. To give a general idea of the VIM, suppose the following general non-linear equation:

$$Lu(x) + Nu(x) = f(x) \quad (7)$$

where $L$ is a linear operator, and $N$ is a non-linear operator. According to the VIM, correction functional can be written in the following form:

$$u_{n+1}(x) = u_n(x) + \int_0^x \lambda(s)[Lu_n(s) + Nu_n(s) - f(s)]ds \quad (8)$$

where $\lambda(s)$ is a general Lagrange multiplier. To make the above correction functional stationary with respect to $u_n$, we have:

$$\delta u_{n+1}(x) = \delta u_n(x) + \delta \int_0^x \lambda(s)[Lu_n(s) + N\tilde{u}_n(s) - f(s)]ds =$$

$$\delta u_n(x) + \int_0^x \lambda(s)\delta[Lu_n(s)]ds = 0 \quad (9)$$

where $\tilde{u}_n$ is considered as a restricted variation i.e. $\delta \tilde{u}_n = 0$. Therefore, we first determine the Lagrange multiplier $\lambda$ that will be identified optimally via integration by parts. The successive approximations $u_{n+1}(x)$ of the solution $u_n(x)$ will be readily obtained upon using the obtained Lagrange multiplier and by using a appropriate function for $u_0(x)$. Therefore by starting from $u_0(x)$, we can obtain the exact solution or an approximate solution of the equation (7) using $u(x) = \lim_{n \rightarrow \infty} u_n(x)$.

## III. SOLUTION OF THE PERONA-MALIK EQUATION BY VIM

Assume the following equation with indicated initial condition:

$$\begin{cases} u_t = \frac{1}{u_x^2 + u_y^2}(u_y^2 u_{xx} - 2u_x u_y u_{xy} + u_x^2 u_{yy}) \\ u(x, y, 0) = \sqrt{x^2 + y^2} - 1 \end{cases} \quad (10)$$

For solving this equation by VIM, correction functional can be written in the following form.

$$u_{n+1}(x,y,t) = u_n(x,y,t)$$
$$+ \int_0^x \lambda(s) \left( \frac{du_n}{ds} - \frac{1}{\left(\frac{\partial u_n}{\partial x}\right)^2 + \left(\frac{\partial u_n}{\partial y}\right)^2} \left( \left(\frac{\partial u_n}{\partial y}\right)^2 \frac{\partial^2 u_n}{\partial x^2} - 2\frac{\partial u_n}{\partial x}\frac{\partial u_n}{\partial y}\frac{\partial^2 u_n}{\partial x \partial y} + \left(\frac{\partial u_n}{\partial x}\right)^2 \frac{\partial^2 u_n}{\partial y^2} \right) \right) ds$$

(11)

To make the above correction functional stationary with respect to $u_n$, we have:

$$\delta u_{n+1}(x,y,t) = \delta u_n(x,y,t)$$
$$+ \delta \int_0^x \lambda(s) \left( \frac{\partial u_n}{\partial x} - \frac{1}{\left(\frac{\partial \widetilde{u_n}}{\partial x}\right)^2 + \left(\frac{\partial \widetilde{u_n}}{\partial y}\right)^2} \left( \left(\frac{\partial \widetilde{u_n}}{\partial y}\right)^2 \frac{\partial^2 \widetilde{u_n}}{\partial x^2} - 2\frac{\partial \widetilde{u_n}}{\partial x}\frac{\partial \widetilde{u_n}}{\partial y}\frac{\partial^2 \widetilde{u_n}}{\partial x \partial y} + \left(\frac{\partial \widetilde{u_n}}{\partial x}\right)^2 \frac{\partial^2 \widetilde{u_n}}{\partial y^2} \right) \right) ds$$
$$= \delta u_n(x,y,t) + \int_0^x \lambda(s) \delta\left(\frac{\partial u_n}{\partial x}\right) ds = 0$$

(12)

From the above relation for any $\delta u_n$, we obtain the Euler-Lagrange equation:

$$\lambda'(s) = 0 \qquad (13)$$

With the following natural boundary condition.

$$\lambda(t) + 1 = 0 \qquad (14)$$

Using equations (8) and (9), Lagrange multiplier can be identified optimally as follows:

$$\lambda(t) = -1 \qquad (15)$$

Substituting the identified Lagrange multiplier into equation (6) results in the following iterative relation:

$$u_{n+1}(x,y,t) = u_n(x,y,t)$$
$$+ \int_0^x \lambda(s) \left( \frac{du_n}{ds} - \frac{1}{\left(\frac{\partial u_n}{\partial x}\right)^2 + \left(\frac{\partial u_n}{\partial y}\right)^2} \left( \left(\frac{\partial u_n}{\partial y}\right)^2 \frac{\partial^2 u_n}{\partial x^2} - 2\frac{\partial u_n}{\partial x}\frac{\partial u_n}{\partial y}\frac{\partial^2 u_n}{\partial x \partial y} + \left(\frac{\partial u_n}{\partial x}\right)^2 \frac{\partial^2 u_n}{\partial y^2} \right) \right) ds$$

(16)

Starting from $u(x,y,0) = \sqrt{x^2+y^2} - 1$, we obtain the following results:

$$u_1 = (x^2 + y^2 + t - \sqrt{x^2+y^2})/(\sqrt{x^2+y^2})$$

$$u_2 = (2x^4 + 2x^2 t - 2x^2\sqrt{x^2+y^2} + 4x^2 y^2 + 2y^4 + 2y^2 t - t^2 - 2\sqrt{x^2+y^2}\, y^2)/(2(x^2+y^2)^{3/2})$$

$$u_3 = (2x^6 + 6x^4 y^2 - 2\sqrt{x^2+y^2}\, x^4 + 2tx^4 + 6y^4 x^2 - 4\sqrt{x^2+y^2}\, x^2 y^2 - x^2 t^2 + 4x^2 y^2 t - 2\sqrt{x^2+y^2}\, y^4 + t^3 + 2y^4 t - t^2 y^2 + 2y^6)/(2(x^2+y^2)^{\frac{5}{2}})$$

$$\vdots$$

(17)

Suppose $u(x,y,t) = u_{13}(x,y,t)$

Finally, After implementation of the obtained solution based on the VIM, we get the results which are shown in the Fig.2. Also, the plot of the VIM solution for t=10 is shown in the Fig.1:

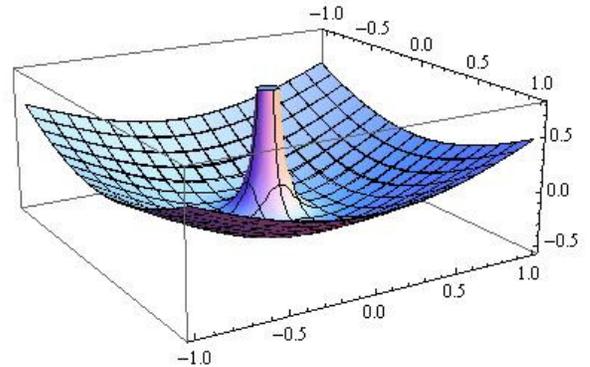

Fig.1. the VIM solution for t=10

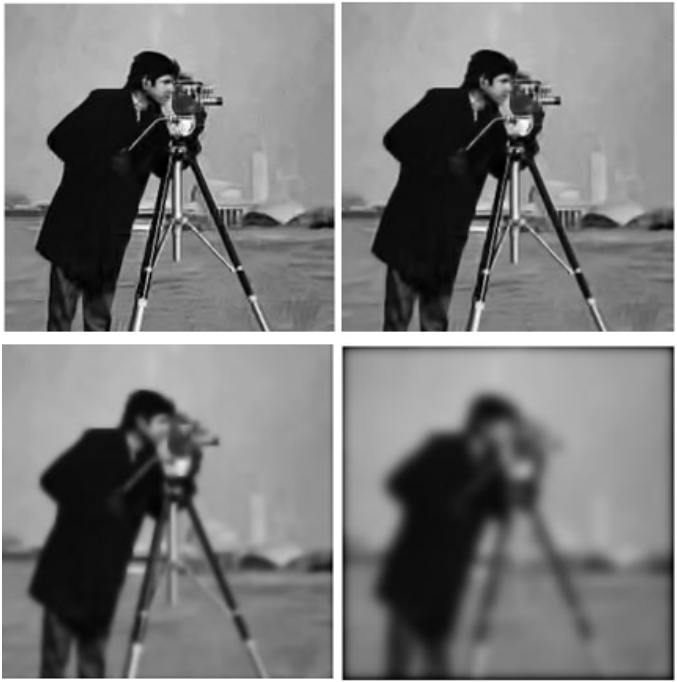

Fig.2. Application of the P-M solution which is obtained via VIM in image restoration.(top, from left to right) a) Cameraman's photo at the original image $t = 0$, b) restored photo after $t = 1$, c) after $t = 10$, d) after $t = 50$.

## IV. CONCLUSION

We could believe that there might be other methods which can solve the problem regularly. But we intend to examine a new numerical method called Radial Basis Functions (RBF's) in this problem. We expect that since the level set method has some valuable properties such as reducing of dimensions, removing singularity and etc. by applying it in association with RBF's we can obtain new optimized algorithm which will lead us to more speedy computations in image restoration. We propose that by elaborating this new combined method, we can assign another application to the level set method and RBF's. It might be used instead of usual methods which people have already applied.

## REFERENCES

[1] E. Nadernejad, H. Koohi, and H. Hassanpour, "PDEs–Based Method for Image Enhancement," *Applied Mathematical Sciences*, Vol. 2, No. 20, pp. 981 – 993, 2008.
[2] Y. You, W. Xu, A. Tannenbaum and M. Kaveh, "Behavioral analysis of anisotropic diffusion in image processing," *IEEE Trans. Image Processing*, vol. 5, No. 11, pp 1539–1553, 1996.
[3] T. Chan, A. Marquina, and P. Mulet, "High–order total variation based image Restoration,' *SIAM J. Sci. Comp*, vol.22, No.2, pp 503–516, 2000.
[4] G. Gilboa, N. Sochen, and Y.Y. Zeevi, "Image Enhancement and Denoising by Complex Diffusion Processes," *IEEE Trans. On Pattern Analysis and Machine Intelligence*, vol. 26, No. 8, 2004.
[5] P. Perona, and J. Malik, "Scale–space and edge detection using anisotropic diffusion," *IEEE Transactions on Pattern Analysis and Machine Intelligence,* vol. 12, No.7, pp. 629–639, 1990.
[6] L. Alvarez, P.L. Lions, and J.M. Morel. "Image selective smoothing and edge detection by nonlinear diffusion.II," *SIAM J. Numer. Anal.*, vol.29,No. 3,pp. 845–866,1992.
[7] z. Kriv´A, "Explicit FV Scheme for the Perona—Malik Equation," *Computational Methods in Applied Mathematics*, vol.5, No.2, pp.170–200, 2005.
[8] L. Rudin, S. Osher, and E. Fatemi, "Nonlinear total variation based noise removal algorithm," *Physica D*, vol. 60, pp. 259–268, 1992.
[9] arXiv:1008.2579v1
[10] G. R. Liu, *Mesh Free Methods: Moving Beyond the Finite Element Method,* CRC Press, Boca Raton, USA, 2002.
[11] P. Kornprobst, R. Deriche, and G. Aubert, "Image sequence restoration: A PDE based coupled method for image restoration and motion segmentation" in Proc. *5th European Conference on Computer Vision*, vol.2, Fribourg, Germany, 1998,pp. 548–562.
[12] L. Rudin and S. Osher, "Total variation based image restoration with free local constraints," in *Proc. IEEE International Conference on Image Processing*, vol. 1, Austin, Texas, 1994, pp. 31–35.
[13] M. Burger, G. Gilboa, S. Osher, and J. Xu, "Nonlinear inverse scale space methods," *Communications in Mathematical Sciences*, vol. 4, pp. 175–208, 2006.
[14] F. Catte, P. L. Lions, J. M. Morel, and T. Coll, "Image selective smoothing and edge detection by nonlinear diffusion," *SIAM J. Numer. Anal.*, vol.29, pp. 182–193, 1992.
[15] G. Aubert and P. Kornprobst, *Mathematical problems in images processing*, Springer–Verlag, Berlin, 2002.
[16] J. H. He, "Variational iteration method- a kind of non-linear analytical technique: some examples," *Int. J. Nonlin. Mech*., vol.34, pp. 699–708, 1999.
[17] J. H. He, "Variational iteration method for autonomous ordinary differential systems," *Appl. Math. Comp.,* vol.114, pp. 115–123, 2000.
[18] J. H. He, "Some asymptotic methods for strongly nonlinear equation," *Int. J. Nod. Phy.,* vol.20, pp. 1144–1199, 2006.
[19] J. H. He and X. H. Wu, "Construction of solitary solutions and compacton-like solution by variational iteration method," *Chos. Soltn. Frcts.,* vol.29, No.1, pp. 108–113, 2006.
[20] M. Inokuti, H. Sekine and T. Mura, "General use of the Lagrange multiplier in non-linear mathematical physics," in Variational Method in the Mechanics of Solids, S. Nemat-Nasser, Ed. Oxford: Pergamon Press, 1978, pp. 156-162.